\documentclass[draft,english]{ourlematema}

\title{Ranks and Singularities \\ of Cubic Surfaces}
\titlemark{Ranks and Singularities of Cubic Surfaces}

\usepackage{graphicx}
\usepackage{array}
\usepackage{color}
\usepackage{pifont}
\newcommand{\cmark}{\ding{51}}%
\newcommand{\xmark}{\ding{55}}%

\def\CC{{\mathbb C}}
\def\RR{{\mathbb R}}
\def \PP{{\mathbb P}}
\def \KK{{\mathbb K}}

\MSC{14J17,
14Q10,
15A69}
\keywords{rank, singularities, cubic surface, Hessian, tensor}

\author{Anna Seigal}
\address{%
Mathematical Institute \\
University of Oxford \\
\email{seigal@maths.ox.ac.uk}
}
\author{Eunice Sukarto}
\address{%
Department of Mathematics \\
University of California, Berkeley \\
\email{eunicesukarto@berkeley.edu}
}

\begin{document}

\maketitle

\begin{abstract}
We explore the connection between the rank of a polynomial and the singularities of its vanishing locus.
We first describe the singularity of generic polynomials of fixed rank. 
We then focus on cubic surfaces.
Cubic surfaces with isolated singularities are known to fall into~22 singularity types. 
We compute the rank of a cubic surface of each singularity type.
This enables us to find the possible singular loci of a cubic surface of fixed rank.
Finally, we study connections to the Hessian discriminant.
We show that a cubic surface with singularities that are not ordinary double points lies on the Hessian discriminant, and that the Hessian discriminant is the closure of the rank six cubic surfaces. 
\end{abstract}

\section{Introduction}

Let $f$ be a homogeneous polynomial of degree $d$ in $n$ variables $x_1, \ldots, x_n$. The rank of $f$ is the smallest number $r$ such that it can be written as a sum
$$ f = l_1^d + \cdots + l_r^d ,$$
where the $l_i$ are linear forms in the variables.
We identify a polynomial $f$ with its vanishing locus.
Recall that a point $p \in \PP^{n-1}$ is a singular point of $f$ if all partial derivatives of~$f$ vanish at $p$. 
All polynomials are assumed to be homogeneous, and we usually require that the coefficients of the linear forms $l_i$ lie in the same field as the coefficients of~$f$.
Polynomials have complex coefficients unless otherwise stated.

Polynomials of degree $d$ in $n$ variables, with coefficients in a field $\KK$ of characteristic zero, are in bijection with symmetric tensors in $(\KK^n)^{\otimes d}$.
This generalises the correspondence between quadratic forms and symmetric matrices to higher degree polynomials.
Symmetric tensors are multi-indexed arrays whose entries are unchanged under permuting the indices. We obtain a polynomial from a tensor $T$ by multiplying the tensor along each index by the vector of indeterminates $x = (x_1, \ldots, x_n)$,
$$ f(x_1, \ldots, x_n) = \sum_{i_1 = 1}^n \cdots \sum_{i_d = 1}^n T_{i_1, \ldots, i_d} x_{i_1} \cdots x_{i_d}.  $$
Conversely, we obtain a symmetric tensor from a polynomial by setting the tensor entry at position $(i_1, \ldots, i_d)$ to be the coefficient of the monomial $x_{i_1} \cdots x_{i_n}$, divided by the multinomial coefficient ${{d \choose j_1, \ldots, j_n}}$ where $x_{i_1} \cdots x_{i_d} = x_1^{j_1} \cdots x_n^{j_n}$. 

The vector space of polynomials of degree $d$ in $n$ variables, considered up to scale, has dimension ${{d + n - 1 \choose d}}-1$.
The generic rank is defined to be the rank of a general polynomial in this space. A formula for the generic rank is
given by the Alexander Hirschowitz theorem, see~\cite{AleHir} and~\cite[Theorem 5.4.1.1]{Lan}. The rank of a polynomial  can exceed the generic rank for degrees $d \geq 3$. 
The polynomials whose rank differs from the generic rank are contained in a hypersurface, by e.g.~\cite[Theorem 1.21]{Sei2}. We will see the connection between this hypersurface and discriminantal loci, for a few examples.

There are a few different notions of rank.
In this paper, we consider the decomposition of a polynomial as a sum of powers of linear forms. This is a decomposition of the corresponding symmetric tensor as a sum of symmetric rank one tensors. 
The rank is the length of a shortest decomposition as a sum of powers of linear forms, the symmetric rank of the corresponding tensor.
A decomposition in which the rank one terms do not need to be symmetric gives the  (non-symmetric) rank.
In the setting of cubic surfaces, the two notions of rank are known to coincide~\cite{Sei}, and they also agree for matrices, though they can differ for a general tensor~\cite{Shi}. 
The border rank is the smallest~$r$ such that a polynomial can be written as a limit of rank $r$ polynomials. For ranks defined over $\RR$, multiple ranks can occur with positive probability when sampling uniformly in the space of polynomials; these are called typical ranks. 

\bigskip

We study the connections between two ways to divide up of the space of polynomials: by rank and by singularities. 
The connection between the rank of a polynomial and its singular points is well-known in the degree two case. The singular locus of a rank $r$ quadratic form is a linear space of codimension $r$. 
We can decompose the quadratic form as $\lambda_1 l_1^2 + \cdots + \lambda_r l_r^2$ where the~$l_i$ are linear forms whose vectors of coefficients $v_i$ are orthogonal, and the $\lambda_i$ are scalars. The singular locus is the orthogonal complement of $\langle v_1, \ldots, v_r \rangle$. 
In particular, a quadratic form is singular if and only if its rank is less than the generic rank, which is $n$. The rank of a quadratic polynomial cannot exceed the generic rank.

We study the connection between rank and singularities for higher degree polynomials.
The first cases are binary cubics and ternary cubics (cubic curves).

\begin{exa}[Binary cubics]
\label{ex:binarycubic}
A summary of the possible ranks and singularities of binary cubics is given in Table~\ref{table:binarycubic}.
 We parametrise binary cubics by 
$f(x_1,x_2) = a_0 x_1^3 + 3 a_1 x_1^2 x_2 + 3 a_2 x_1 x_2^2 + a_3 x_2^3$
for $(a_0 :a_1: a_2: a_3) \in \PP^3$.
The singular binary cubics lie on the discriminant hypersurface
$$ 3 a_1^2 a_2^2 + 6 a_0 a_1 a_2 a_3 - 4 a_1^3 a_3 - 4 a_2^3 a_0 - a_0^2 a_3^2 .$$
The discriminant is the set of binary cubics not of generic rank.
\begin{table}[ht]
\centering
\begin{tabular}{c|c|c}
 Rank & Normal form & Singularities \\ \hline
1 & $x_1^3$ & 1 \\ 
2 & $x_1^3 + x_2^3$ & 0 \\
3 & $x_1^2 x_2$ & 1 \\
\end{tabular}
\caption{Ranks and singularities of binary cubics}
\label{table:binarycubic}
\end{table}
\end{exa}

\begin{exa}[Cubic curves]
\label{ex:cubiccurve}
The rank and singularities of cubic curves
are given in Table~\ref{table:ternarycubic}. The ranks can be found in~\cite[Section 10.4]{Lan} or~\cite[Section 8]{LanTei}, and the singularity types in~\cite{BruceWall}. See also~\cite[\S96]{Segre}.
\begin{table}[ht]
\centering
\begin{tabular}{c|c|c}
Normal form & Rank & Singularities \\ \hline
$x_1^3$ & 1 & $\PP^1$ \\ 
$x_1^3 + x_2^3$ & 2 & $D_4$ \\
$x_1^2 x_2$ & 3 & $\PP^1$ \\
$x_1^3 + x_2^3 + x_3^3$ & 3 & 0 \\
$x_1^3 + x_2^3 + x_3^3 - 3 \rho x_1x_2x_3$, $\rho^3 \neq 1, 0, -8$ & 4 & 0 \\
$x_2^2x_3 - x_1^3 - x_1^2x_3$ & 4 & $A_1$ \\
$x_2^2x_3 - x_1^3$ & 4 & $A_2$ \\
$x_1(x_1^2 + x_2x_3)$ & 4 & $2A_1$ \\
$x_1x_2x_3$ & 4 & $3A_1$ \\
$x_2(x_1^2 + x_2x_3)$ & 5 & $A_3$ \\
\end{tabular}
\caption{Ranks and singularities of cubic curves}
\label{table:ternarycubic}
\end{table}
The number of singular points that can occur for each rank is summarised in Table~\ref{table:ternaryversus}.
The discriminant of singular cubic curves is a degree 12 hypersurface, see e.g.~\cite[Chapter 1, 4.15]{GKZ}. The discriminant contains all cubic curves whose rank exceeds the generic rank. We observe that only cubic curves of generic rank four can have multiple isolated singular points. 
\begin{table}[ht]
\centering
\begin{tabular}{c|l|l|l|l|l}
\multicolumn{1}{c|}{} & \multicolumn{5}{c}{No. singular points}                    \\ \hline
Rank & 0   & 1   & 2 & 3  & $\infty$ \\ \hline
1   & \xmark  & \xmark & \xmark & \xmark & \cmark   \\
2    &  \xmark & \cmark  &  \xmark & \xmark & \xmark \\
3   & \cmark  & \xmark & \xmark & \xmark & \cmark   \\
4   & \cmark & \cmark & \cmark& \cmark & \xmark \\
5   & \xmark & \cmark & \xmark & \xmark & \xmark   \\
\end{tabular}
\caption{Rank versus number of singular points for cubic curves}
\label{table:ternaryversus}
\end{table}
\end{exa}

In this article, we extend these examples to cubic surfaces. We also give a partial generalisation to general polynomials. 

\bigskip

One connection between the rank of a polynomial and the singular points of its vanishing locus is given in~\cite{LanTei}.
Given a polynomial $f$ of degree $d$ in $n$ variables, its flattening $f_{k,d-k}$ is the matrix of size $n^k \times n^{d-k}$ obtained from the symmetric tensor corresponding to $f$ by combining the first $k$ indices to index the rows, and taking the remaining $d-k$ indices to index the columns. 
Let $\Sigma_s(f)$ denote the set of points that vanish to order $s$ on $f$. The authors show that the rank of $f$ is bounded strictly from below by $\text{rank} f_{s,d-s} + \dim \Sigma_s (f)$.
This generalises the fact that the rank is bounded from below by the rank of a flattening.
In the context of cubic polynomials, this says that the rank is strictly greater than the flattening rank plus the dimension of the space of singular points. In particular, a singular cubic minimally written using $n$ variables has rank at least $n+1$. In~\cite{DimSti}, the authors study the types of singularities that occur using the hyperplane arrangement given by the linear forms in a decomposition into rank one terms.

Ranks and singularities have also been studied in the context of real rank boundaries, the hypersurfaces that separate loci of different typical real ranks.
The discriminant is conjectured to be a defining polynomial of certain real rank boundaries~\cite{MicMoo,MicMoo2}.
Discriminants of binary cubics describe the real rank two polynomials~\cite{SeiStu}. For binary forms, one real rank occurs generically on each connected component of the complement of the discriminant~\cite{Ble,ComOtt}. 

\bigskip

The rest of this paper is organised as follows.
We consider how the rank of a polynomial relates to its membership in the discriminant hypersurface in Section~\ref{sec:fixedrank}.
We focus on cubic surfaces 
in Section~\ref{sec:isolated}. We find the rank of a cubic surface of each singularity type. We use this to find the combinations of rank and singular points that can occur for cubic surfaces, and
we characterise the possible singular loci of general polynomials whose rank is at most the number of variables. 
We also find the possible dimensions and degrees of infinite singular loci of cubic surfaces. We show that the highest-dimensional orbit of cubic surfaces with infinitely many singular points has rank six.
We study connections to the Hessian discriminant in Section~\ref{section:hd}. 
We show that a singular polynomial is in the Hessian discriminant unless all singularities are ordinary double points, and we show that the Hessian discriminant is the closure of the rank six cubic surfaces. This implies that a general singular cubic surface has rank five.

\section{Generic polynomials of fixed rank}
\label{sec:fixedrank}

A general rank $r$ tensor is 
a linear combination of $r$ powers of general linear forms, provided $r$ is at most the generic rank. It is harder to find a polynomial of rank exceeding the generic rank: in most cases the existence of such a polynomial is not known~\cite{BucTei}.
A generic polynomial of rank $r$ also has border rank $r$ for $r$ at most the generic rank, but not for ranks exceeding the generic rank. 

The rank one polynomials are powers of linear forms $f =~l^d$. Considered up to scale, they parametrise the Veronese variety $\nu_d(\PP^{n-1})$. 
A weighted sum of $r$ general points on the Veronese variety parametrises the $r$th secant variety. 
The $r$th secant variety is the closure of the rank $r$ tensors, considered up to scale, provided the rank is at most the generic rank. 
The generic rank is the smallest $r$ such that the $r$th secant variety of the Veronese fills the space of polynomials. 

In this section, we describe the singularity of a generic polynomial of fixed rank.
We give the containment relations between the discriminant and secant varieties of the Veronese.
We show that a general rank $r$ polynomial in $n$ variables is singular if $r <n$ and non-singular for ranks from $n$ to the generic rank.
We give a condition for a rank $n+1$ polynomial in $n$ variables to be singular.
The results in this section are expository, though they may not have been written down elsewhere. 

\begin{prop}
\label{prop:contains}
Let $f$ be a homogeneous polynomial of degree~$d$ in $n$ variables, and let $f = l_1^d + \cdots + l_r^d$ be a minimal length decomposition as a sum of powers of linear forms. Then the singular locus of $f$ contains a linear space of codimension~$r$. If the coefficients of the $l_i$ are linearly independent vectors, there are no other singular points. 
\begin{proof}
The partial derivatives of $f$ 
vanish if $l_i(p) = 0$ or, equivalently, $\langle v_i, p \rangle = 0$ where $v_i$ is the vector of coefficients of $l_i$. If the $v_i$ are linearly independent, the singular locus is the orthogonal complement of $\langle v_1, \ldots, v_r \rangle$, since we can apply a general linear group transformation so that $f = \sum_{i = 1}^r x_i^d$. 
\end{proof}
\end{prop}

For a generic polynomial in $n$ variables of rank at most $n$, the rank one terms in a decomposition are linearly independent and Proposition~\ref{prop:contains} gives all singular points. When the rank exceeds the number of variables, Proposition~\ref{prop:contains} does not contribute any singular points, but singular points can still arise.

\begin{prop}
Let
\begin{equation}
f = a_1 x_1^d + \cdots + a_n x_n^d + a_{n+1} (-x_1 - \cdots - x_n)^d 
\label{eqn:frankn+1}
\end{equation}
be a polynomial in $n \geq 2$ variables $x_i$ of degree $d \geq 2$ with coefficients $a_j \in \CC^*$. 
Then $f$ is non-singular unless
\begin{equation}
\label{eqn:rankn+1}
\lambda_{n+1}^{d-1} = {\left( - \sum_{j = 1}^n \lambda_j \right)}^{d-1},
\end{equation}
where $\lambda_j^{d-1} = \frac{a_1}{a_j}$, and $f$ is singular at exactly the points $(\lambda_1 : \cdots : \lambda_n)$ for which Equation~\eqref{eqn:rankn+1} holds.
\begin{proof}
At a singular point, $a_ix_i^{d-1}=a_jx_j^{d-1}$ for all $1\leq i,j\leq n$.
Since not all coordinates $x_i$ vanish, these conditions imply that no coordinates vanish, and we set $x_1 = 1$ to obtain $x_j^{d-1} = \frac{a_1}{a_j}$. Setting $\lambda_j$ to be the value of $x_j$ at the singular point recovers the condition in Equation~\eqref{eqn:rankn+1}.
\end{proof}
\label{prop:rankn+1}
\end{prop}

A test for the singularity of a polynomial is given by the vanishing of the discriminant, but an expression for the discriminant is not known in general~\cite[Chapter 13]{GKZ}.
The singularity of a polynomial of the form in Equation~\eqref{eqn:frankn+1} can be checked by computing
Equation~\eqref{eqn:rankn+1} for all choices of $\lambda_j$. There are $d-1$ choices for $\lambda_j$, for $1 \leq j \leq n$, hence there are $(d-1)^n$ conditions to check. The existence of a real singular point can be checked via $2^n$ conditions if $d$ is odd, and a single condition if $d$ is even, since there are two real $(d-1)$st roots of $\frac{a_1}{a_j}$ for odd $d$, and a single real root for even $d$.

\begin{thm}
\label{thm:nleqrleqrg}
Let $f$ be a general polynomial of rank $r$, where $r$ is at least the number of variables and at most the generic rank. Then $f$ is non-singular.
\begin{proof}
The polynomial $f$, up to scale, is a general point on the $r$th secant variety of the Veronese. 
The secant variety is irreducible, hence~$f$ is non-singular provided the secant variety is not contained in the discriminant. 
It suffices to observe that the non-singular polynomial $x_1^d + \cdots + x_n^d$ is in the secant variety.
\end{proof}
\end{thm}

\begin{cor}
A general polynomial in $n$ variables of border rank $r$ is singular if $r < n$ and non-singular if $r \geq n$.
\begin{proof}
 A general polynomial of border rank $r$ is a general point on the $r$th secant variety of the Veronese. 
 For $r<n$ the result follows from Proposition~\ref{prop:contains}. 
 For $r \geq n$ the result follows from Theorem~\ref{thm:nleqrleqrg}.
\end{proof}
\end{cor}

\section{Ranks and singularity types}
\label{sec:isolated}

In this section we first consider the ranks of cubic surfaces with isolated singularities.
Cubic surfaces with isolated singularities fall into 22 singularity types~\cite{BruceWall,Sch}.
Our main result is Table~\ref{table:isolatedcubic}, which lists the ranks of a cubic surface of each type.
There is a single normal form for each type of singularity, in all but one case.
There are two normal forms with a $D_4$ singularity, see~\cite[Lemma 4]{BruceWall}. The two normal forms can be distinguished by their rank.
We summarise the possible number of singular points for a cubic surface of each rank in Table~\ref{table:maintab}.
We give a partial generalisation to general polynomials in Table~\ref{table:tabn}.
We study the possible dimensions and degrees of infinite singular loci of cubic surfaces in Table~\ref{table:degreeanddim}.
We show that a general cubic surface with infinitely many singular points has rank six.

\begin{tiny}
\begin{table}[htbp]
\centering
\addtolength{\tabcolsep}{-5pt}
\begin{tabular}{c|c|c|c}
No. & Sing.  & Normal form & Rank \\ \hline \hline
1  & $A_1$             & $\begin{matrix} x_4(x_2^2-x_1x_3) + x_2x_1^2-(1+\rho_1+\rho_2+\rho_3)x_1x_2^2 \\ + (\rho_1+\rho_2\rho_3)x_1x_2x_3 +(1+\rho_1)(\rho_2+\rho_3)x_2^3 \\ - ((1+\rho_1)\rho_2\rho_3+\rho_1(\rho_2+\rho_3))x_2^2x_3 + \rho_1\rho_2\rho_3x_2x_3^2 \\ \rho_i\in \mathbb{C}\backslash\{0,1\} \text{ pairwise different}\end{matrix}$ & 
$\begin{matrix} 5 \text{ for } \rho \\ \text{ generic e.g.} \\
(-1,2,-2) \end{matrix}
$
\\
\hline
2  & $2A_1$            & $\begin{matrix} x_4(x_2^2-x_1x_3) + x_2(x_1-(1+\rho_1)x_2+\rho_1x_3)(x_1-\rho_2x_2) \\ \rho_i\in \mathbb{C}\backslash\{0,1\} \text{ not equal} \end{matrix}$ & $\begin{matrix} 5 \text{ for } \rho \\ \text{ generic e.g.} \\
(-1,2) \end{matrix}
$ \\ \hline
3  & $3A_1$            & $\begin{matrix} x_4(x_2^2-x_1x_3) + x_2^2(x_1-(1+\rho)x_2+\rho x_3) \\ \rho \in \mathbb{C}\backslash\{0,1\}\end{matrix}$ & $\begin{matrix} 5 \text{ if} \\ 
\rho^2 \neq \frac{14}{9}\! \rho \! - \! 1 \end{matrix} $
\\ \hline
4  & $4A_1$            & $x_4(x_2^2-x_1x_3)+x_2^2(x_1-2x_2+x_3)$ & 5 \\ \hline
5  & $A_1A_2$          & $\begin{matrix} x_4(x_2^2-x_1x_3)+x_1x_2(x_1-(1+\rho)x_2+\rho x_3) \\ \rho \in \mathbb{C}\backslash \{0,1\} \end{matrix}$ & $\begin{matrix} 6 \text{ if } \rho \text{ is} \\ \frac12(1 \pm \sqrt{3} i) \end{matrix}$ \\ \hline
6  & $2A_1A_2$         & $x_4(x_2^2-x_1x_3)+x_2^2(x_1-x_2)$ & 6 \\ \hline
7  & $A_12A_2$         & $x_4(x_2^2-x_1x_3)+x_2^3$ & 6 \\ \hline
8  & $A_1A_3$          & $x_4(x_2^2-x_1x_3)+x_1^2x_2-x_1x_2^2$ & 6 \\ \hline
9  & $2A_1A_3$         & $x_4(x_2^2-x_1x_3)+x_1x_2^2$ & 6 \\
\hline
10 & $A_1A_4$          & $x_4(x_2^2-x_1x_3)+x_1^2x_2$ & 6 \\
\hline
11 & $A_1A_5$          & $x_4(x_2^2-x_1x_3)+x_1^3$ & 6 \\
\hline
12 & $A_2$             & $\begin{matrix} x_4x_1x_2-x_3(x_1^2+x_2^2+x_3^2+\rho_1 x_1x_3+\rho_2x_2x_3) \\ \rho_1, \rho_2 \in \mathbb{C}\backslash \{\pm 2\} \end{matrix}$ & $6 \text{ if } \rho_i = 0$ \\
\hline
13 & $2A_2$            & $x_4x_1x_2-x_3(x_1^2+x_3^2+\rho x_1x_3), \rho \in \mathbb{C}\backslash \{\pm 2 \}$ & $6 \text{ if } \rho ^2\neq 3$ \\
\hline
14 & $3A_2$            & $x_4x_1x_2-x_3^3$ & 5 \\
\hline
15 & $A_3$             & $x_4x_1x_2+x_1(x_1^2-x_3^2)+\rho x_2(x_2^2-x_3^2), \rho \neq 0$ & 6 \\
\hline
16 & $A_4$             & $x_4x_1x_2+x_1^2x_3+x_2(x_2^2-x_3^2)$ & 6 \\
\hline
17 & $A_5$             & $x_4x_1x_2+x_1^3+x_2(x_2^2-x_3^2)$ & 6 \\
\hline
18 & $D_4^{I}$            & $x_4x_1^2+x_2^3+x_3^3+x_1x_2x_3$ & 6 \\
\hline
19 & $D_4^{II}$           & $x_4x_1^2+x_2^3+x_3^3$ & 5 \\
\hline
20 & $D_5$             & $x_4x_1^2+x_1x_3^2+x_2^2x_3$ & 6 \\
\hline
21 & $E_6$             & $x_4 x_1^2+x_1x_3^2+x_2^3$ & 6 \\
\hline
22 & $\tilde{E_6}$ & $x_1^3+x_2^3+x_3^3-3\rho x_1x_2x_3, \, \rho^3\neq 1$ & $\begin{matrix} 3, \rho^3 = 0,-8 \\ 4, \text{otherwise} \end{matrix}$ \\
\end{tabular}
\caption{Ranks of cubic surfaces with isolated singularities}
\label{table:isolatedcubic}
\end{table}
\end{tiny}

The ranks of singular cubic surfaces of each singularity type are computed as follows.
By~\cite{LanTei}, the rank is at least five if all four variables are used in the normal form, i.e. if the corresponding tensor has full flattening rank. By~\cite[Theorem 2.10]{Sei}, the rank of a cubic surface with finitely many singular points is at most six. Hence, in most cases, we seek to distinguish between ranks five and six. 

We obtain lower bounds using the substitution method, see~\cite[\S5.3.1]{Lan}. Cubic surfaces correspond to symmetric $4 \times 4 \times 4$ tensors. If the $4 \times 4$ slices are linearly independent, and no multiple of three of the slices can be subtracted from the fourth to give a rank two matrix, then the rank of the tensor is at least six. We also obtain lower bounds by subtracting other combinations of slices.

We upper bound the rank by showing the existence of a decomposition, using the structure of the Hessian.
The Hessian of a generic cubic surface
has 10 ordinary double points, where the rank of the $4 \times 4$ Hessian matrix drops from three to two.
These 10 points are the triple intersections of the linear forms in the unique rank five decomposition:
if the cubic surface has decomposition $f = l_1^3 + \cdots + l_5^3$, then the 10 points are the triple intersections $l_i = l_j = l_k = 0$. 
This is Sylvester's pentahedral theorem, see~\cite[\S{}84]{Segre}. 
Sylvester's pentahedron is the five planes $l_i = 0$, as well as their double and triple intersections. 
The decomposition can therefore be obtained from the singular points of the Hessian. We use this idea to find decompositions of some singular cubic surfaces.

The normal forms in Table~\ref{table:isolatedcubic} are from~\cite{Sch}, except with fewer parameters, since we consider normal forms with respect to the general linear group instead of the special linear group, because rank is unchanged under general linear transformation. 
We now describe the 22 cases. 
For cases 1, 2, 3, and 13, we find the rank of a generic cubic surface with that singularity type. For cases 5 and 12 we give parameter values for which we computed the rank. In all other cases, we compute the rank of every cubic surface with that singularity type.

\begin{enumerate}
\item[1.] 
We equate the normal form when $\rho = (-1,2,-2)$ to a sum of five cubes of linear forms.
A Macaulay2~\cite{M2} computation shows that the rank is five:
the system of equations is in 14 parameters and has codimension 14, indicating a finite number of rank five decompositions. Specialising the linear forms in the decomposition to be compatible with the singular points of the Hessian allows the computation to terminate. The rank is five for generic $\rho$ by the results in Section~\ref{section:hd}.

\begin{footnotesize}
\begin{verbatim}
R = QQ[a1,a2,a4,b1,b2,b4,c2,c4,d1,d2,d4,e1,e2,e4][x1,x2,x3,x4];
f = x4*(x2^2 - x1*x3) + x2*(x1-x3)*(x1-4*x3);
g1 = (a1*x1 + a2*x2 - 2*a1*x3 + a4*x4)^3;
g2 = (b1*x1 + b2*x2 - 2*b1*x3 + b4*x4)^3;
g3 = (c2*x2 + c4*x4)^3;
g4 = (d1*x1 + d2*x2 + 2*d1*x3 + d4*x4)^3;
g5 = (e1*x1 + e2*x2 + 2*e1*x3 + e4*x4)^3;
F = f - (g1 + g2 + g3 + g4 + g5);
(M,C) = coefficients F; I = minors(1,C);
codim I  
\end{verbatim}
\end{footnotesize}

\item[2.] We equate the normal form with parameters $\rho$ to a sum of five cubes of linear forms where, for three linear forms in the decomposition, the coefficient of $x_3$ is $-\rho_1 \rho_2$ times the coefficient of $x_4$. We check that this has a solution when, e.g., $\rho = (-1,2)$ using Macaulay2. The rank being five for generic $\rho$ follows from Section~\ref{section:hd}.

\item[3.] 
There exists a rank five decomposition
$$ \begin{matrix} (24 \rho) f & = &
 (x_1 - \rho x_3 + x_4)^3 + (x_1 + \rho x_3 - x_4)^3 + (-x_1 + \rho x_3 - x_4)^3 \\ & & + \lambda (x_1 + \alpha x_2 + \rho x_3 + x_4)^3 \\
 & & -(1 + \lambda) (x_1 - (3(\rho + 1) - \alpha) x_2 + \rho x_3 + x_4)^3. \end{matrix}$$
 The decomposition was obtained by finding linear forms whose triple intersections give singular points of the Hessian of $f$.
The variables $\alpha$ and $\lambda$ are solutions to
$\alpha^2+3(\rho + 1)\alpha + 8 \rho$ and $\lambda (2 \alpha + 3 \rho + 3 ) +\alpha + 3 \rho+3$, 
which have solutions provided $2 \alpha + 3\rho + 3 \neq 0$, i.e. $9 \rho^2 - 14 \rho + 9 \neq 0$.

\item[4.] We have the rank five decomposition
$$\begin{matrix} 24f & = & (-x_1 + x_3 + x_4)^3 + (x_1 - x_3 + x_4)^3 +
(x_1 + x_3 - x_4)^3 \\ & & - 2(x_1 - 2x_2 + x_3 + x_4)^3 + (x_1 - 4x_2 + x_3 + x_4)^3 .\end{matrix} $$

\item[5.] Denote the four slices of the tensor by $M_i$. We construct the $4 \times 4 \times 2$ tensor with slices $M_1 = M_1 + s_1 M_3 + s_2 M_4$ and $M_2 = M_2 + t_1 M_3 + t_2 M_4$. 
It has $4 \times 8$ flattening
\begin{footnotesize}
$$ \begin{bmatrix} 
0 & \rho s_1 + 2 & -s_2 & -s_1 & 2 & \rho t_1-2 \rho-2 & \rho-t_2 & -t_1 \\
      \rho s_1+2 & -2 \rho+2 s_2-2 & \rho & 0 & \rho t_1-2 \rho-2 & 2 t_2 & 0 & 2 \\ -s_2 & \rho &
      0 & -1 & \rho-t_2 & 0 & 0 & 0 \\ -s_1 & 0 & -1 & 0 & -t_1 & 2 & 0 & 0 \\
\end{bmatrix} $$
\end{footnotesize}
When $\rho^2 - \rho + 1 = 0$, there do not exist values of $s_1,s_2,t_1,t_2$ such that the $4 \times 4$ minors of the flattening vanish. Hence the rank of the $4 \times 4 \times 2$ tensor is at least four, and the rank of the original tensor is six.

\item[6.] Similarly to 5, we construct the $2 \times 4 \times 4$ tensor with slices $M_1 - s_1 M_2 - t_1 M_3$ and $M_4 - s_2 M_2 - t_2 M_3$. It has $4 \times 8$ flattening
$$ \begin{bmatrix} 0 & -2 s_1 & 0 & t_1 & 0 & -2 s_1 & -1 & t_2 \\
-2 s_1 & 2 + 6 s_1 & 0 & -2 s_1 & -2 s_2 & 2 + 6 s_2 & 0 & -2 s_2 \\
0 & 0 & 0 & -1 & -1 & 0 & 0 & 0 \\
t_1 & -2 s_1 & -1 & 0 & t_2 & -2 s_2 & 0 & 0 \\
\end{bmatrix} .$$
The minor on columns $1,3,4,7$ is $2s_1$, and the minor on columns $2,3,4,7$ is $-2(1 + 3 s_1)$. They cannot both vanish, hence the matrix has rank four and the original tensor has rank six. 

\item[7.] The slices are linearly independent. No multiple of the first, second, and third slices can be subtracted from the fourth to give a rank two matrix, so the rank is at least six. 

\item[8.] Assume that the cubic surface has rank $r$. The $4 \times 4 \times 2$ tensor with slices $M_1 + s_1 M_3 + s_2 M_4$ and $M_2 + t_1 M_3 + t_2 M_4$ has rank at most $r-2$. There exists a multiple of the $4 \times 2$ matrix $T_{4jk}$ that can be subtracted from the slices $T_{ijk}$ for $1 \leq i \leq 3$ to give a $3 \times 4 \times 2$ tensor of rank at most $r-3$. Its $4 \times 6$ flattening is
$$ \begin{bmatrix} * & * & * & * & * & * \\ 2 & -2 & * & * & 0 & 2 \\ 
* & * & 0 & 0 & -1 & 0 \\
* & * & * & 2 & * & 0 \end{bmatrix}, $$
where the $*$ represent unknown entries. This matrix does not have a rank two completion, hence $r-3 \geq 3$ and the original rank is six. 

    \item[9.] Subtract an arbitrary multiple of the second and third slices from the first and fourth slices to give a $4 \times 4 \times 2$ with $4 \times 8$ flattening 
    $$ \begin{bmatrix} 0 & * & 0 & * & 0 & * & -1 & * \\
    * & 2 & 0 & * & * & 2 & 0 & * \\
    0 & 0 & 0 & -1 & -1 & 0 & 0 & 0 \\
    * & * & -1 & 0 & * & * & 0 & 0 \\
    \end{bmatrix},$$
    This cannot be completed to a rank three matrix, hence $r-2 \geq 4$ and the original rank is six. 
    
    \item[10.] No multiple of the first, second, and third slices can be subtracted from the fourth to give a matrix with minor $M_{123}$, on rows $1,2,3$ and columns $1,2,3$, vanishing. So the rank is at least six.
    \item[11.] Same as 10.
    \item[12.] First we subtract multiple $s_i$ of the fourth slice from the $i$th slice for $i~\in~\{ 1, 2, 3 \}$ to give a $4 \times 4 \times 3$ tensor. Then we subtract multiple $t_i$ of the fourth slice in another direction to give the $4 \times 3 \times 3$ tensor
    $$ \left[ \begin{array}{ccc||ccc||ccc}
0 & t_1 & -2 & s_1 & s_2 + t_2 & s_3 & -2 & t_3 & -2 \rho_1 \\
s_1 + t_1 & s_2 & s_3 & t_2 & 0 & -2 & t_3 & -2 & -2 \rho_2 \\
-2 & 0 & -2 \rho_1 & 0 & -2 & -2 \rho_2 & -2 \rho_1 & -2 \rho_2 & -6 \\
0 & 1 & 0 & 1 & 0 & 0 & 0 & 0 & 0 
\end{array} \right] ,$$
where $||$ separates the slices of the tensor. Finally we consider the matrix obtained by subtracting an arbitrary multiple of the first and second slices from the third slice. The $2 \times 2$ minors of this $4 \times 3$ matrix cannot vanish if, e.g., $\rho_1 = \rho_2 = 0$. Hence the matrix has rank at least two, and the original tensor has rank six. Other parameters for which this works include $(\rho_1,\rho_2) = (0,\pm \sqrt{3})$ and $(\pm \sqrt{3}, 0)$. 
    \item[13.] Subtracting a multiple of the second, third, and fourth slices from the first slice and taking the ideal of $3 \times 3$ minors gives the condition $\rho^2 = 3$ so, if this does not hold, the rank is at least six.
    \item[14.] The first monomial has rank four and the second has rank one, so the rank is at most five. 
    \item[15.] Subtracting an arbitrary multiple of the third and fourth slices from the first and second slices gives a tensor with $4 \times 8$ flattening matrix
    $$ \begin{bmatrix} 
    6 & * & * & 0 & 0 & * & * & 1 \\
    * & 0 & * & 1 & * & 6 \rho & * & 0 \\
    * & * & -2 & 0 & * & * & -2 \rho & 0 \\
    0 & 1 & 0 & 0 & 1 & 0 & 0 & 0 
    \end{bmatrix}.$$
    It does not have a rank three completion, hence the original rank is six.
    \item[16.] No multiple of the first, third, and fourth slices can be subtracted from the second to give a vanishing minor $M_{134}$, so the rank is at least six. 
    \item[17.] Same as 16.
    \item[18.] No multiple of the second, third, and fourth slices can be subtracted from the first to give a vanishing minor on rows $1,3,4$ and columns $1,2,4$, so the rank is at least six.
    \item[19.] The normal form, written as a tensor, has only five non-vanishing entries, so the rank is at most five. 
    \item[20.] No multiple of the second, third, and fourth slices can be subtracted from the first to give a vanishing minors $M_{134}$, so the rank is at least six. 
    \item[21.] Same as 20. 
    \item[22.] See the fourth and fifth normal forms for cubic curves in Table~\ref{table:ternarycubic}. When $\rho^3 = 1$ the curve is singular, and the surface has infinitely many singular points. When $\rho^3 = 0$ or $-8$ the cubic curve is non-singular with singular Hessian, so it has the fourth normal form. For all other $\rho$ the curve and its Hessian are both non-singular and it has the fifth normal form. 
\end{enumerate}

The computation of the rank of a cubic surface of each singularity type allows us to tabulate the possible combinations of rank and singularities. 
 We see that, as for cubic curves, only the cubic surfaces with generic rank, in this case rank five, can have the maximal number of isolated singular points.

\begin{thm}
\label{thm:main}
The possible combinations of rank and singular points, for a cubic surface, are as in Table~\ref{table:maintab}. 

\begin{table}[htbp]
\centering
\begin{tabular}{llllll}
\multicolumn{1}{l|}{}     & \multicolumn{5}{l}{No. singular points}                                                                  \\ \hline
\multicolumn{1}{l|}{Rank} & \multicolumn{1}{l|}{0}   & \multicolumn{1}{l|}{1}   & \multicolumn{1}{l|}{2}   & \multicolumn{1}{l|}{3}   & 4   \\ \hline
\multicolumn{1}{c|}{1}    & \multicolumn{1}{l|}{\xmark}  & \multicolumn{1}{l|}{\xmark}  & \multicolumn{1}{l|}{\xmark}  & \multicolumn{1}{l|}{\xmark}  &  \xmark  \\
\multicolumn{1}{c|}{2}    & \multicolumn{1}{l|}{\xmark}  & \multicolumn{1}{l|}{\xmark}  & \multicolumn{1}{l|}{\xmark}  & \multicolumn{1}{l|}{\xmark}  &  \xmark  \\
\multicolumn{1}{c|}{3}    & \multicolumn{1}{l|}{\xmark}  & \multicolumn{1}{l|}{\cmark} & \multicolumn{1}{l|}{\xmark}  & \multicolumn{1}{l|}{\xmark}  &  \xmark  \\
\multicolumn{1}{c|}{4}    & \multicolumn{1}{l|}{\cmark} & \multicolumn{1}{l|}{\cmark} & \multicolumn{1}{l|}{\xmark}  & \multicolumn{1}{l|}{\xmark}  &  \xmark  \\
\multicolumn{1}{c|}{5}    & \multicolumn{1}{l|}{\cmark} & \multicolumn{1}{l|}{\cmark} & \multicolumn{1}{l|}{\cmark} & \multicolumn{1}{l|}{\cmark} & \cmark \\
\multicolumn{1}{c|}{6}    & \multicolumn{1}{l|}{\cmark} & \multicolumn{1}{l|}{\cmark} & \multicolumn{1}{l|}{\cmark} & \multicolumn{1}{l|}{\cmark}   & \xmark   \\
\multicolumn{1}{c|}{7}    & \multicolumn{1}{l|}{\xmark}  & \multicolumn{1}{l|}{\xmark}  & \multicolumn{1}{l|}{\xmark}  & \multicolumn{1}{l|}{\xmark}  &   \xmark
\end{tabular}
\caption{The possible combinations of rank and number of singular points}
\label{table:maintab}
\end{table}

\begin{proof}
There cannot be a finite number of singular points for ranks one and two by Proposition~\ref{prop:contains}, and
for rank seven by~\cite[Theorem 2.10]{Sei}.
For rank three, there is one singular point when the three linear forms are linearly independent. If the linear forms are dependent, there are an infinite number of singular points. For rank four, if the four linear forms are independent, the cubic surface is non-singular. If they are linearly dependent, the surface is a cone over a cubic curve, say in variables $(x_1, x_2, x_3)$. If the curve is singular, then the surface has an infinite number of singular points. If the cubic curve is non-singular, the only singular point of the surface is $(0:0:0:1)$. 

A general cubic surface of rank five is non-singular.
Table~\ref{table:isolatedcubic} shows the existence of rank five cubic surfaces with one, two, three and four singular points. 
Non-singular cubic surfaces of rank six can be found in~\cite{Segre}. Singularity types $A_4$, $A_1A_4$, and $A_12A_2$ have rank six. There does not exist a cubic surface of rank six with exactly four singular points, because $4A_1$ is the only singularity type with four singularities, and its rank is five.
\end{proof}
\end{thm}

We give a partial generalisation to polynomials in $n$ variables.

\begin{thm}
\label{thm:generaln}
The possible singular points of a homogeneous polynomial of degree $d \geq 3$
in $n \geq 4$ variables,
with rank at most $n$, are as shown in Table~\ref{table:tabn}.

\begin{table}[ht]
\centering
\begin{tabular}{ccccc}
\multicolumn{1}{c|}{}     & \multicolumn{4}{c}{Number of singular points}                                                                  \\ \hline
\multicolumn{1}{c|}{Rank} & \multicolumn{1}{c|}{$0$}   & \multicolumn{1}{c|}{$1$}   & 
\multicolumn{1}{c|}{$1 < k < \infty$} & {$\infty$}  \\ \hline
\multicolumn{1}{c|}{$\leq n-2$}  & \multicolumn{1}{c|}{\xmark}  & \multicolumn{1}{c|}{\xmark}  &  \multicolumn{1}{c|}{\xmark} & \cmark  \\
\multicolumn{1}{c|}{$n-1$}    & \multicolumn{1}{c|}{\xmark}  & \multicolumn{1}{c|}{\cmark}  &  \multicolumn{1}{c|}{\xmark} & \cmark  \\
\multicolumn{1}{c|}{$n$}    & \multicolumn{1}{c|}{\cmark}  & \multicolumn{1}{c|}{\cmark} &  \multicolumn{1}{c|}{\xmark} & \cmark  \\
\end{tabular}
\caption{The possible singular points of a polynomial in $n$ variables of rank $\leq n$}
\label{table:tabn}
\end{table}

\begin{proof}
If the linear forms appearing in a decomposition span a space of dimension $k \leq n-2$, the singular locus is positive dimensional by Proposition~\ref{prop:contains}. 
Hence a polynomial of rank at most $n-2$ cannot have finitely many singular points.
The singular locus cannot contain isolated singular points, because no conditions are imposed on the orthogonal complement of the linear forms, so any singular point lies in a $n-k-1$ dimensional linear space of singular points.

A polynomial of rank $n-1$ with finitely many singular points is projectively equivalent to $x_1^d + \dots + x_{n-1}^d$, with one singular point.
In the rank~$n$ case, if the span of the linear forms has dimension $n-1$, then $f$ is a cone over a polynomial~$g$ in $n-1$ variables of rank $n$. By Theorem~\ref{thm:nleqrleqrg}, there exists a non-singular such~$g$ provided $n$ does not exceed the generic rank, for which the hypothesis $n \geq 4$ is sufficient. Then $f$ has one singular point. If the span of the linear forms has dimension $n$, the polynomial is projectively equivalent to $x_1^d + \dots + x_n^d$, with no singular points.

It remains to construct rank $n-1$ and rank $n$ polynomials with infinitely many singular points.
Consider
$$ x_1^d + \cdots + x_{n-2}^d - a_{n-1} (x_1 + \cdots + x_{n-2})^d .$$
It has rank $n-1$ provided $a_{n-1} \neq 0$.
Setting $a_{n-1} = \frac{1}{(n-2)^{d-1}}$ gives a singular point whenever $x_1 = \cdots = x_{n-2} = 1$, for any values of $x_{n-1}$ and $x_n$. The analogous example in one more variable gives a rank $n$ example with infinitely many singular points. 
\end{proof}
\end{thm}

We now consider cubic surfaces with infinitely many singular points. We characterise the possible dimensions and degrees of the singular locus.

\begin{thm}
\label{thm:infinitesing}
The possible dimensions and degrees of infinite singular loci of a cubic surface are as shown in Table~\ref{table:degreeanddim}. 

\begin{table}[ht]
\centering
\addtolength{\tabcolsep}{-3pt}
\begin{tabular}{c|c|c|c|c|c}
Rank & dim 1, deg 1 & dim 1, deg 2 & dim 1, deg 3 & dim 2, deg 1 & other \\ \hline
1    & \xmark           & \xmark   &   \xmark     & \cmark          & \xmark           \\
2    & \cmark          & \xmark   &  \xmark      & \xmark           & \xmark           \\
3    & \xmark          & \xmark   &   \xmark     & \cmark          & \xmark           \\
4    & \cmark          & \cmark   &  \cmark      & \xmark           & \xmark           \\
5    & \cmark          & \xmark   &   \xmark     & \xmark           & \xmark           \\
6    & \cmark          & \cmark  &   \xmark     & \xmark           & \xmark           \\
7    & \cmark          & \cmark  &  \xmark      & \xmark           & \xmark          
\end{tabular}
\caption{Dimensions and degrees of the infinite singular loci}
\label{table:degreeanddim}
\end{table}

\begin{proof}
The result for ranks one and two follows from the proof of Proposition~\ref{prop:contains}.
For rank three, if the three linear forms are linearly independent then the cubic is singular at a single point. Hence they must span a two-dimensional subspace, and the cubic has normal form $x_1^2 x_2$, with singular points the plane in $\PP^3$ with equation $x_1 = 0$. For rank four, the four linear forms appearing in a decomposition can span at most a three-dimensional space, and we have a cone over a cubic curve. The singular locus of a rank four cubic curve ranges from zero points to three points, see Table~\ref{table:ternarycubic}, hence the infinite singular loci of the cone range from one line to three lines.

Any rank five cubic, minimally expressed using four variables, can be written
$$ f = a_1 x_1^3 + a_2 x_2^3 + a_3 x_3^3 + a_4 x_4^3 + a_5 (-x_1 - \cdots - x_m)^3,  $$
where $a_i \in \CC^*$. Note that this differs from the general form in Equation~\eqref{eqn:frankn+1} because we may have $m < n$. The singular points are the singular points of
$$ f = a_1 x_1^3 + \cdots + a_m x_m^3 + a_5 (-x_1 - \cdots - x_m)^3,  $$
the finitely many points given in Proposition~\ref{prop:rankn+1}, with the added condition $x_k = 0$ for $k > m$. Hence there are finitely many singular points overall. If the span is three dimensional, we have a cone over the cubic curve with normal form $x_2(x_1^2 + x_2 x_3)$, which has a single singular point, hence the cubic surface has a line of singular points. 

In the rank six case, there are three cases to consider, with normal forms $x_1(x_1^2 + x_2^2 + x_3^2 + x_4^2)$, $x_1(x_2^2 + x_3^2 + x_4^2)$, and $x_1 x_2^2 + x_3 x_4^2$, see~\cite{Segre,Sei}. In the first two cases the singular locus is a curve of degree two, and in the third case it is a line. In the rank seven case, the possible normal forms are $x_1 (x_1 x_2 + x_3^2 + x_4^2)$ and $x_1^2 x_2 + x_1 x_3 x_4 + x_3^3$, see~\cite{CGV,Segre}. In the first case, the singular locus is a curve of degree two, and in the second case it is a line. 
\end{proof}
\end{thm}

\begin{cor}
\label{cor:infinite}
On the highest-dimensional orbit of cubic surfaces with infinitely many singular points, a cubic surface has rank six.
\begin{proof}
There are finitely many orbits of cubic surfaces with infinitely many singular points, see the proof of Theorem~\ref{thm:infinitesing}. We seek the orbit whose dimension is largest. The stabiliser of the rank six cubic surface $x_1 x_2^2 + x_3 x_4^2$ has lower dimension than the other normal forms, hence the orbit has highest dimension, by~\cite[Proposition 21.4.3]{TauYu}.
\end{proof}
\end{cor}

We seek the implicit description of the rank six cubic surfaces, in order to obtain the rank of a general singular cubic surface.
The variety of polynomials whose rank exceeds the generic rank is known in only a few cases, see~\cite{BucTei} and Examples~\ref{ex:binarycubic} and~\ref{ex:cubiccurve}. 
We see that the rank six cubic surfaces parametrise a discriminantal locus known as the Hessian discriminant.

\section{The Hessian discriminant} 
\label{section:hd}

The Hessian matrix of a polynomial is its matrix of second order partial derivatives. The Hessian variety is the hypersurface defined by the determinant of the Hessian matrix. 
For a cubic polynomial, the Hessian is a symmetric matrix of linear forms.

Consider the hypersurface defined by the determinant of a $n \times n$ symmetric matrix of ${{n + 1 \choose 2}}$ indeterminate entries. The singular points occur when the partial derivatives of the determinant vanish, that is, when the $(n-1) \times (n-1)$ minors vanish. This is a variety of codimension 3 
and degree ${{ n + 1 \choose 3}}$
in the space of symmetric matrices~\cite{Harris}.
For example, for a $4 \times 4$ symmetric matrix the variety has degree 10. 

The Hessian variety of a generic cubic surface
has 10 ordinary double points, where the rank of the $4 \times 4$ Hessian matrix drops from three to two. 
The Hessian discriminant is the locus of cubic surfaces whose Hessian matrix does not intersect the rank two matrices at 10 reduced points.
It was introduced in~\cite{Sei}, and described in more detail in~\cite{DinuSey} and~\cite[Chapter 5]{Sei2}. 
By~\cite{Stu17}, the Hessian discriminant is the Hurwitz form to the rank two symmetric matrices, specialised to the linear space of symmetric matrices that occur as the Hessian of a cubic surface. 

In this section we study how the rank and singularity type of a cubic surface interface with its membership in the Hessian discriminant. 

\begin{prop}
\label{prop:HD}
A cubic surface with singularities that are not ordinary double points lies on the Hessian discriminant. \begin{proof}
For a cubic surface $f$, we denote by $J_f$ the radical of the ideal of $3 \times 3$ minors of the Hessian matrix. 
For the orbits whose normal forms do not involve parameters, including all those with infinite singular locus, we check directly that $J_f$ does not consist of exactly 10 points. The cases involving parameters are $A_1A_2$, $A_2$, $2A_2$, $A_3$, and $\tilde{E_6}$. For these cases, we compute the degree of $J_f$ for a general choice of parameters via a multidegree computation~\cite{MilStu}. We consider the polynomial ring in which the variables $x_i$ have degree $(1,0)$ and the parameters $\rho_j$ have degree $(0,1)$. The degree for general parameters is the part of the multidegree only involving the degree of variables $x_i$. In all four cases the ideal $J_f$ has degree at most five in the variables $x_i$,
hence $J_f$ does not consist of 10 points and the singularity types lie in the Hessian discriminant.
\end{proof}
\end{prop}

Cubic surfaces with no singularities or only ordinary double points are those that are stable, in the sense of geometric invariant theory, under the action of $SL_4$, see~\cite[Theorem 3.6]{Rei}.
The intersection of the Hessian discriminant and the discriminant contains the cubic surfaces of singularity types 5 to 22 in Table~\ref{table:isolatedcubic}, by Proposition~\ref{prop:HD}. The null cone is contained in this intersection. It consists of those with singularities not of type~$A_1$ or $A_2$, types 8 to 11 and 15 to 22, see~\cite{Rei}.

The cubic surfaces with $k$ ordinary double points have normal forms expressed in terms of $4-k$ parameters, see Table~\ref{table:isolatedcubic}. Proposition~\ref{prop:HD} shows that, for almost all parameters,
the polynomial does not lie on the Hessian discriminant.
However, for some special parameters,
cubic surfaces with only ordinary double points can lie on the Hessian discriminant.

\begin{exa}
We find the cubic surfaces of singularity type $3A_1$ that lie on the Hessian discriminant. 
Singularity type $3A_1$ has normal form 
$f = x_4(x_2^2-x_1x_3) + x_2^2(x_1-(1+\rho)x_2+\rho x_3)$,
where $\rho \in \mathbb{C}\backslash\{0,1\}$. The ideal of $3 \times 3$ minors of the Hessian of $f$ decomposes into 
four points
$$ ( -\rho : 0 : 1 : 0), \quad 
( 0 : 1 : 0 : 0 ), \quad
( 0 : 0 : 1 : -\rho), \quad 
(1 : 0 : 0 : -1),$$
and three ideals of degree two, parametrised by
$$ (\rho : x : 1 : 0) , \quad (0 : x : 1 : \rho), \quad (1: y : 0 : 1) ,$$
where $4 x^2 - 3(\rho + 1) x + 2 \rho = 0$ and $4 \rho y^2 - 3(\rho + 1)y + 2 = 0$. 
The ideal does not consist of 10 distinct points when the discriminant $9 \rho^2 - 14 \rho + 9$ vanishes. The singularity type for these parameters is still $3A_1$, because for all three singular points the affine chart with that singular point at the origin gives a  non-homogeneous polynomial whose degree two part has rank three~\cite{BruceWall}.
\end{exa}

We now discuss how rank relates to the Hessian discriminant.
In~\cite{DinuSey}, the authors show that the defining equation of the Hessian discriminant is the cube of the degree 40 Salmon invariant of a cubic surface, and that it is parametrised by general rank six non-singular cubic surfaces, which have normal form
$x_1^3 + x_2^3 + x_3^3 + x_4^2 (\lambda_1 x_1 + \lambda_2 x_2 + \lambda_3 x_3 + \lambda_4 x_4)$ for $\lambda_i \in \CC$.
It is the boundary of the real rank five cubic surfaces, described in~\cite{MicMoo}.
In~\cite[\S6.6]{DarGee}, the authors show that the vanishing of the degree 40 Salmon invariant parametrizes the cubic surfaces of rank greater than five. See also~\cite[\S9.4.5]{Dolgachev}. Together, the results of~\cite{DarGee,DinuSey} rule out the existence of cubic surfaces of rank five whose decomposition is not obtained from the singular points of its Hessian, and implies the following. 

\begin{thm}
The Zariski closure of the set of rank six cubic surfaces is the Hessian discriminant.
\end{thm}

In the previous section, we constructed rank five decompositions of some singular cubics using the 10 points where the Hessian has rank two.
We now extend this to show that a rank five decomposition can always be built when the Hessian drops to rank two at exactly 10 reduced points, a direct argument for the results described above from~\cite{DarGee,Dolgachev}.
Note that the Hessian is singular at these 10 points, as well as at the singular points of $f$.

\begin{thm}
\label{thm:sylfail}
A cubic surface not on the Hessian discriminant has rank five. 
The Hessian discriminant is the locus of cubic surfaces for which Sylvester's pentahedral theorem fails. 
\begin{proof}
Let $f$ be a cubic surface not on the Hessian discriminant.
The border rank of $f$ is five, since
the fourth secant variety of the Veronese lies in the Hessian discriminant: the Hessian of $x_1^3 + x_2^3 + x_3^3 + x_4^3$ is rank two along six lines where two coordinates vanish.

Write $f$ as a limit of rank five cubic surfaces $f = \lim_{t \to 0} f(t)$, where $f(t) = l_1(t)^3 + \cdots + l_5(t)^3$.
Following~\cite[\S{}89]{Segre}, the $f(t)$ can be chosen such that the planes $l_i(t) = 0$, and their double and triple intersections, have well-defined limit.
Let $l_i = \lim_{t \to 0} l_i(t)$.
The Hessian of $f$ has rank two at the triple intersections of the $l_i$, since the rank of a matrix cannot increase in a limit.
The Hessian of $f$ has rank two at exactly 10 reduced points, which are therefore the distinct triple intersections of the $l_i$. Any four~$l_i$ are linearly independent, since otherwise multiple triple  intersections would coincide. In particular no two planes $l_i=0$ can be the same.

The rank of a polynomial exceeds the border rank when it lies on a limit of planes through $r$ points on the Veronese, but not on a plane through $r$ points.  In order for the rank to exceed the border rank, the points $l_i^3$ on the Veronese must be linearly dependent~\cite[Section 10]{LanTei}. 
Without loss of generality, we have $l_1^3 = a_2 l_2^3 + \cdots + a_k l_k^3$, where $2 \leq k \leq 5$. Since $l_2, \ldots, l_k$ are linearly independent, the only scalars $a_2, \ldots, a_k$ that give a linear power is if all the $a_i$ vanish except one. Then two limiting planes are the same, a contradiction.

The locus of points where Sylvester's pentahedral theorem fails contains the Hessian discriminant~\cite{Segre}. The above shows that the locus where the theorem fails is contained in the Hessian discriminant, hence the two are equal.
\end{proof}
\end{thm}

Theorem~\ref{thm:sylfail} is shown for smooth cubic surfaces in~\cite[\S89]{Segre}.
It is a discriminantal description of the locus of cubic surfaces of non-generic rank. It would be interesting to obtain such discriminants for other spaces of polynomials.

We conclude with the following result, which follows from the fact that a general point on the discriminant is not on the Hessian discriminant.

\begin{cor}
A generic singular cubic surface has rank five.
\end{cor}

{\bf Acknowledgements.}
The authors thank Jennifer Sixt Pinney and Bernd Sturmfels for instigating our first meeting. 
ES thanks Holly Mandel and AS thanks David Eisenbud and Giorgio Ottaviani for helpful discussions.

\vspace{-5ex}
\nocite{*}

\end{document}